\documentclass[11pt]{article}
\usepackage{mathrsfs}
\usepackage{amsfonts}
\usepackage{amsmath,amssymb}
\newcommand{\epsln}{\varepsilon}

\newcommand{\rmd}{\mbox{\rm d}}
\newcommand{\bfa}{{\mbox{\boldmath $a$}}}
\newcommand{\bfb}{{\mbox{\boldmath $b$}}}
\newcommand{\bfc}{{\mbox{\boldmath $c$}}}
\newcommand{\bff}{{\mbox{\boldmath $f$}}}

\newcommand{\bfn}{{\mbox{\boldmath $n$}}}
\newcommand{\bfu}{{\mbox{\boldmath $u$}}}
\newcommand{\bfv}{{\mbox{\boldmath $v$}}}

\newcommand{\bfr}{{\mbox{\boldmath $r$}}}
\newcommand{\bfx}{{\mbox{\boldmath $x$}}}

\newcommand{\bfB}{{\mbox{\boldmath $B$}}}
\newcommand{\bfR}{{\mbox{\boldmath $R$}}}

\newcommand{\bfJ}{{\mbox{\boldmath $J$}}}
\newcommand{\bfK}{{\mbox{\boldmath $K$}}}

\newcommand{\bfM}{{\mbox{\boldmath $M$}}}

\newcommand{\bfE}{{\mbox{\boldmath $E$}}}
\newcommand{\bfF}{{\mbox{\boldmath $F$}}}

\def\bgeq{\begin{equation}}
\def\edeq{\end{equation}}
\def\bgar{\begin{array}}
\def\edar{\end{array}}

\title{Modification to Maxwell's Equations}
\author{{Shangbin Cui}\\
}
\date{}

\begin{document}
\maketitle

\begin{abstract}
  A fixed electric charge is an electric current relative to a moving magnetic field, so that it
  is subjected to the force of the moving magnetic field. This means that not only time-varying
  magnetic field produces electric field, but moving magnetic field produces electric field as
  well. Maxwell neglected this fact in deriving his equations for the description of dynamical
  behavior of electromagnetic field so that the two equations (A) $\nabla\cdot\bfE=\epsln_0^{-1}\rho$
  and (B) $\partial_t\bfE-\epsln_0^{-1}\mu_0^{-1}\nabla\times\bfB=-\epsln_0^{-1}\bfJ$ are incorrect.
  In this paper we modify the equation (A) into $\nabla\cdot(\bfE+\overline{\bfu}\times\bfB)=
  \epsln_0^{-1}\rho$, where $\overline{\bfu}$ denotes the mean velocity of the charges in the
  electric current, and the equation (B) is correspondingly modified. The modified equations are
  invariant under Galilean transformation. As a byproduct of this work, we see that Einstein's
  theory of special relativity is wrong.
\medskip

\textbf{Keywords}:\ Modification; Maxwell's Equations; electromagnetic field; Galilean transformation;
  invariance.
\medskip

\end{abstract}

\section{Introduction}

\hskip 2em
  Maxwell's equations describing the dynamical behavior of electromagnetic field induced by electric
  charges and their current are one of the most important components of theoretical physics and also
  the whole science. In the case that the space is a vacuum, these equations read as follows (cf., e.g.,
  \cite{Jackson}):
$$
  \partial_t\bfB+\nabla\times\bfE=\bf0,
\eqno{(1.1)}
$$
$$
  \partial_t\bfE-\epsln_0^{-1}\mu_0^{-1}\nabla\times\bfB=-\epsln_0^{-1}\bfJ,
\eqno{(1.2)}
$$
$$
  \nabla\cdot\bfE=\epsln_0^{-1}\rho,
\eqno{(1.3)}
$$
$$
  \nabla\cdot\bfB=0,
\eqno{(1.4)}
$$
$$
  \partial_t\rho+\nabla\cdot\bfJ=0.
\eqno{(1.5)}
$$
  Here $\bfE$ and $\bfB$ denote the electric field strength and the magnetic flux density, respectively,
  $\rho$ and $\bfJ$ are the density of electric charges and the strength of the current of these charges,
  respectively, by which the field $(\bfE,\bfB)$ is induced, and $\epsln_0$ and $\mu_0$ are the
  dielectric constant in vacuum and the magnetic permeability constant of the vacuum, respectively. The
  last equation is the mathematical expression of the law of conservation of charges, so that it is obvious.
  The other four equations were derived by James Clark Maxwell from the Coulomb's Law describing the
  behavior of electrostatic field and its generalization to non-electrostatic field, the Biot-Savart Law
  describing the behavior of stable magnetic field, the Faraday's Law of electromagnetic induction
  which states that a time-varying magnetic field induces an electric field, and the hypothesis that also
  a time-varying electric field induces a magnetic field. To be more precise, let us apply the Helmholtz
  decomposition theorem for vector field to write
$$
  \bfE=\bfE_g+\bfE_c \qquad \mbox{and} \qquad \bfB=\bfB_g+\bfB_c,
$$
  where the subscripts $g$ and $c$ indicate the gradient parts and the curl parts, respectively, of the
  corresponding vector fields on the left sides of these equations. The equations (1.1) and (1.3) show that
  $\bfE_c$ and $\bfE_g$ are respectively induced by time-variation of magnetic field and the distribution
  of charges in the space, the equation (1.2) shows that electric current and time-variation of electric
  field jointly induce $\bfB_c$, and the equation (1.4) shows that $\bfB_g\equiv\bf0$.

  Maxwell proposed the above equations in 1865 \cite{Maxwell}. These equations elevated the study of
  electromagnetism from the experimental stage onto a comprehensive theoretical stage.
  In particular, since Heinrich Rudolf Hertz confirmed in 1887 the existence of electromagnetic waves
  predicted by Maxwell from the above equations, the study of electromagnetics has changed the world.
  Nowadays we humans live in a world filled with various man-made electromagnetic waves and our daily
  lives are completely inseparable from electromagnetic waves. Such an outcome was entirely beyond
  anticipation in the era before Maxwell.

  A byproduct of Maxwell's equations was the emergence of Albert Einstein's theory of special relativity.
  In 1887, Albert A. Michelson and Edward W. Morley conducted an experiment that demonstrated non-existence
  of the ether, a medium for propagating electromagnetic waves hypothesized by Maxwell,
  which implies that an absolutely stationary frame of reference for validity of these equations does not
  exist. Nonexistence of an absolutely stationary frame of reference implies that all inertial frame of
  references should be equal in position for the description of the Maxwell's equations, so that these
  equations should be invariant under Galilean transformation of the frame of references. Unfortunately,
  not all these equations are invariant under Galilean transformation. This leads to a paradox: Since
  non-invariance of these equations under Galilean transformation implies that observations of the same
  electromagnetic phenomena obtained by two observers in relative motion are different, whose observation
  should we adopt? In 1892 Hendrik Lorentz resolved this issue by introducing the idea to set different
  scales for length and time in different reference frames (A similar idea was proposed by another
  physicist George Francis Fitzgerald in 1889). More precisely, Lorentz found that the equations
  (1.1)--(1.5) were invariant under the so-called Lorentz transformation. On the basis of the idea of Lorentz,
  in 1905 Einstein proposed the theory of special relativity which has now been one of the most important
  theories in modern physics.

  One may easily verify that the equations (1.1), (1.4) and (1.5) are invariant under the Galilean
  transformation (see Section 3) and the ones that break invariance are equations (1.2) and (1.3). Actually,
  these last two equations are incorrect because they ignore the effect of the motion of magnetic
  field on electric field. Indeed, a moving magnetic field induces an additional electric field and this
  additional electric field has additional reaction on the production of magnetic field, whereas in these
  two equations this issue is completely neglected. More precisely, a fixed electric charge is an electric
  current relative to a moving magnetic field, so that it is subjected to the force of the moving magnetic
  field. Hence, not only time-varying magnetic field induces electric field, but moving magnetic field
  induces electric field as well. The purpose of this paper is to modify the equations (1.2) and (1.3) by
  taking this fact into consideration. After modification, the new Maxwell's equations have the following
  expressions:
$$
  \partial_t\bfB+\nabla\times\bfE=\bf0,
\eqno{(1.6)}
$$
$$
  \partial_t\bfE-\overline{\bfu}\times(\nabla\times\bfE)+\dot{\overline{\bfu}}\times\bfB
  -\nabla\times[\overline{\bfu}\times(\bfE+\overline{\bfu}\times\bfB)]
  -\epsln_0^{-1}\mu_0^{-1}\nabla\times\bfB=-\epsln_0^{-1}\bfJ,
\eqno{(1.7)}
$$
$$
  \nabla\cdot(\bfE+\overline{\bfu}\times\bfB)=\epsln_0^{-1}\rho,
\eqno{(1.8)}
$$
$$
  \nabla\cdot\bfB=0,
\eqno{(1.9)}
$$
$$
  \partial_t\rho+\nabla\cdot\bfJ=0,
\eqno{(1.10)}
$$
  where $\overline{\bfu}=\overline{\bfu}(t)$ denotes the mean velocity of charges in the electric current.
  All the above five equations are invariant under the Galilean transformation of the frame of references.

  In summary, in deriving the equations (1.1)--(1.5) Maxwell used the following ideas:

  (1)\ An electric charge induces an electric field (Coulomb's Law).

  (2)\ An electric current induces a magnetic field (Biot--Savart's Law).

  (3)\ Time-variation of a magnetic field induces electric field (Faraday's Law).

  (4)\ Time-variation of an electric field induces magnetic field (Maxwell's Hypothesis).

\noindent
  In this paper we modify these equations by taking the following fact into consideration:

  (5)\ Motion of a magnetic field also induces electric field.

\noindent
  We then obtain the equations (1.6)--(1.10). The equations (1.1)--(1.5) are invariant under
  Lorentz transformation of the frame of references, but fail to keep invariant under Galilean
  transformation. Unlike this, the equations (1.6)--(1.10) are invariant under Galilean  transformation.

  From the expressions of the equations (1.6)--(1.10) we see that if an electromagnetic system (i.e.,
  the system of $\rho$ and $\bfJ$) has the nice property that the mean velocity $\overline{\bfu}=
  \overline{\bfu}(t)$ is constant, i.e., $\overline{\bfu}(t)\equiv\bfu_0$ for some constant vector
  $\bfu_0$, then in the frame of reference in which $\bfu_0=\bf0$, the equations (1.6)--(1.10) are
  respectively the same as the equations (1.1)--(1.5), so that there exists a special frame of reference
  in which the old Maxwell theory of electrodynamics works. Besides, in our real world that electromagnetic 
  waves are extensively employed, the values of $|\bfu|$ are neglectfully small relative to the light 
  speed $c$, so that the equations (1.6)--(1.10) can be very well approximated by the equations (1.1)--(1.5), 
  application of the old Maxwell theory does not produce any negative effect.

  The rest part of this paper is arranged as follows: In the next section we give the derivation of
  the equations (1.7) and (1.8). In Section 3 we show that the equations (1.6)--(1.10) are invariant
  under Galilean transformation. In the last section we make analysis to a concrete example to show that
  Lorentz transformation fails to apply to electromagnetic issues so that the theory of special
  relativity established on the basis of it does not work. More discussions on invalidation of this
  theory and related problems in theoretical physics will be made in a subsequent paper \cite{Cui} and
  this paper is actually cut from it.

\section{Modification to Maxwell's Equations}

\hskip 2em
  In this section we derive the equations (1.7) and (1.8).

  We use the notation $\bfK$ to denote the Euclidean solid space, i.e., $\bfK={\bfR}^3$ when a Cartesian
  solid coordinate system is established. We assume that the whole space $\bfK$ is a vacuum. Let $P$ denote
  a point variable in $\bfK$, i.e., $P=(x,y,z)$ when a Cartesian solid coordinate system is established. Let
  $\rho=\rho(P,t)$ be the density of electric charges distributed in $\bfK$, and $\bfJ=\bfJ(P,t)$ be the
  strength of electric current in $\bfK$. By definition,
$$
  \bfJ(P,t)=\rho(P,t)\bfu(P,t),
\eqno{(2.1)}
$$
  where $\bfu(P,t)$ is the velocity of the electric charge at the point $P$ at time $t$.

  We first derive the equation (1.8). Let us consider the electromagnetic field induced by a charge with
  quantity $q>0$. If this charge is fixed at the point $O$, then it only induces an electric field in the
  space which by Coulomb's Law has the following expression:
$$
  \bfE(P)=\dfrac{q\bfr}{4\pi\epsln_0r^3}, \qquad P\in\bfK\backslash\{O\},
\eqno{(2.2)}
$$
  where $\bfr=\overrightarrow{OP}$ and $r=|\bfr|$. If this charge moves in a velocity $\bfu$, then apart
  from the above electric field which we re-denote as $\bfE_1$, an additional electric field which we
  denote as $\bfE_2$ is induced. To see this let us consider a test charge with quantity $Q$ placed
  at an arbitrary point $P$. We know that the motion of the charge $q$ induces a magnetic field $\bfB$
  which by the Biot-Savart Law has the expression
$$
  \bfB(P)=\dfrac{\mu_0 q}{4\pi r^3}\bfu\times\bfr, \qquad P\in\bfK\backslash\{O\}.
\eqno{(2.3)}
$$
  Since the charge $q$ is moving, it follows that this magnetic field is moving and, consequently, the
  fixed charge $Q$ at the point $P$ induces a electric current $-Q\bfu\delta(M-P)$ relative to this
  magnetic field (note that we are only considering the value of this current at time $t=0$), where $M$
  denotes an arbitrary point on the line $L$ of direction $\bfu$ across the point $P$, and $\delta$
  denotes the one-dimensional Dirac delta function, so that the charge $Q$ undergoes a force $\bff$:
$$
  \bff=-\int_{L} Q\bfu\delta(M-P)\rmd s(M)\times\bfB(P)=-Q\bfu\times\bfB(P)
  =-\dfrac{\mu_0 Qq}{4\pi r^3}\bfu\times(\bfu\times\bfr),
\eqno{(2.4)}
$$
  where $\rmd s(M)$ denotes the arc length element on $L$. Note that this result can also be obtained
  directly from applying the Lorentz force formula (i.e., the first equality and the first expression 
  following it can be omitted). Hence
$$
  \bfE_2(P)=-\bfu\times\bfB(P)=-\dfrac{\mu_0 q}{4\pi r^3}\bfu\times(\bfu\times\bfr)
  =\dfrac{\mu_0 q}{4\pi r^3}[|\bfu|^2\bfr-(\bfu\cdot\bfr)\bfu], \quad P\in\bfK\backslash\{O\}.
\eqno{(2.5)}
$$
  It follows that when the charge $q$ moves in a velocity $\bfu$, the electric field induced by it in the
  space is
\begin{eqnarray*}
  \qquad\qquad\qquad\;\; \bfE(P) &\,=\,&\bfE_1(P)+\bfE_2(P)
  =\dfrac{q\bfr}{4\pi\epsln_0r^3}-\bfu\times\bfB(P)
\\
  &\,=\,&\dfrac{q\bfr}{4\pi\epsln_0r^3}
  +\dfrac{\mu_0 q}{4\pi r^3}[|\bfu|^2\bfr-(\bfu\cdot\bfr)\bfu].
\\
  &\,=\,&\dfrac{q\bfr}{4\pi\epsln_0r^3}(1+\beta^2)-\dfrac{\mu_0 q}{4\pi r^3}(\bfu\cdot\bfr)\bfu,
  \qquad P\in\bfK\backslash\{O\},
   \qquad\qquad\qquad\quad\;\;  (2.6)
\end{eqnarray*}
  where $\beta=c^{-1}|\bfu|$ and $c$ is the velocity of light in vacuum (recall that $c=1/\sqrt{\mu_0\epsln_0}$).
  From this expression of $\bfE$ we see that for the electromagnetic field $(\bfE,\bfB)$ induced by a moving
  charge $q$ in a velocity $\bfu$ there holds
$$
  \bfE(P)+\bfu\times\bfB(P)=\dfrac{q\bfr}{4\pi\epsln_0r^3}, \qquad  P\in\bfK\backslash\{O\},
\eqno{(2.7)}
$$
  so that
$$
  \nabla\cdot[\bfE(P)+\bfu\times\bfB(P)]=\epsln_0^{-1}q\delta(P), \qquad  P\in\bfK.
$$
  Here $\delta$ denotes the three-dimensional Dirac delta function. For the general case that charges are
  distributed in a bounded domain $\Omega$ with density $\rho=\rho(M)$ and velocity $\bfu=\bfu(M)$
  ($M\in\Omega$), we denote by $(\bfE_M(P),\bfB_M(P))\rmd V(M)$ the electromagnetic field induced by the
  charge $\rho(M)\rmd V(M)$ and corresponding current at the point $M\in\Omega$, where $\rmd V(M)$ denotes
  the volume element at the point $M$. By noticing that
$$
  \bfE(P)=\int_{\Omega}\bfE_M(P)\rmd V(M), \quad\quad
  \int_{\Omega}\!\epsln_0^{-1}\rho(M)\delta(P\!-\!M)\rmd V(M)=\epsln_0^{-1}\rho(P)
$$
  and using the standard method of integration we obtain the following equation:
$$
  \nabla\cdot[\bfE(P)\!+\!\int_{\Omega}\!\bfu(M)\!\times\!\bfB_M(P)\rmd V(M)]
  =\epsln_0^{-1}\rho(P) \quad  \mbox{for}\;\, P\in\bfK.
\eqno{(2.8)}
$$
  Note that this equation can also be directly derived from the equation (2.7) by using the divergence
  theorem as in traditional derivation of the equation (1.3) from the Coulomb's Law (cf., e.g., the
  deduction of (1.13) on Pages 27$\sim$29 of \cite{Jackson}). More precisely, for each point $M\in\Omega$
  we use the equation (2.7) to get
$$
  \bfE_M(P)\rmd V(M)+\bfu(M)\times\bfB_M(P)\rmd V(M)=\dfrac{\rho(M)\overrightarrow{MP}}{4\pi\epsln_0|\overrightarrow{MP}|^3}\rmd V(M),
  \qquad  P\in\bfK\backslash\{M\},
$$
  This implies that for every closed smooth surface $S$ we have
$$
  \int_{S}[\bfE_M(P)\rmd V(M)+\bfu(M)\!\times\!\bfB_M(P)\rmd V(M)]\cdot\bfn(P)\rmd S(P)
  \!=\!\left\{
\begin{array}{ll}
  \epsln_0^{-1}\!\rho(M)\rmd V(M), & \;\; \mbox{if}\; M\in D\\
   \;\;0, & \;\; \mbox{if}\; M\not\in D,
\end{array}
\right.
$$
  where $\bfn(P)$ denotes the unit outward normal of the surface $S$ at the point $P$, $\rmd S(P)$ is the
  surface measure element on $S$, and $D$ is the domain enclosed by $S$. Integrating both sides of the above
  equation with respect to $M$ on $\Omega$, we get
$$
  \int_{S}\Big[\bfE(P)+\Big(\int_{\Omega}\bfu(M)\!\times\!\bfB_M(P)\rmd V(M)\Big)\Big]
  \cdot\bfn(P)\rmd S(P)=\int_{D}\epsln_0^{-1}\!\rho(M)\rmd V(M).
$$
  By using the divergence theorem and the arbitrariness of the closed surface $S$, we obtain (2.8).
  We now define the magnetic field $\bfB$ as follows:  First let
$$
  \overline{\bfu}=\dfrac{1}{|\Omega|}\int_{\Omega}\bfu(M)\rmd V(M).
\eqno{(2.9)}
$$
  Next we define $\bfB$ to be the field such that the following relation holds true:
$$
  \overline{\bfu}\times\bfB(P)=\int_{\Omega}\bfu(M)\times\bfB_M(P)\rmd V(M),  \qquad \forall P\in\bfK.
\eqno{(2.10)}
$$
  Then the equation (2.8) reduces into the equation (1.8).

  We note that in the above deduction we omitted time variable $t$ for simplicity of notations. In what
  follows we recover notations containing the variable $t$, so that
$$
  \bfE=\bfE(P,t), \quad \bfB=\bfB(P,t), \quad \rho=\rho(P,t), \quad \bfu=\bfu(P,t), \quad
  \bfJ=\bfJ(P,t), \quad \overline{\bfu}=\overline{\bfu}(t).
$$
  By differentiating both sides of (1.8) in $t$ and using (1.6) and (1.10), we get
$$
  \nabla\cdot[\partial_t\bfE-\overline{\bfu}\times(\nabla\times\bfE)+\dot{\overline{\bfu}}\times\bfB
  +\epsln_0^{-1}\bfJ]=0.
$$
  Hence, there exists a vector field $\bfM$ such that
$$
  \partial_t\bfE-\overline{\bfu}\times(\nabla\times\bfE)+\dot{\overline{\bfu}}\times\bfB
  +\epsln_0^{-1}\bfJ=\nabla\times\bfM.
$$
  Considering the equation (1.2) and the requirement that the modified equation should be invariant
  under Galilean transformation (see the next section), we let the vector field $\bfM$ have the
  following expression:
$$
  \bfM=\overline{\bfu}\times(\bfE+\overline{\bfu}\times\bfB)+\epsln_0^{-1}\mu_0^{-1}\bfB.
$$
  In this way we obtain the equation (1.7).

  Since the equations (1.1) and (1.4) are compatible, from the derivation of the equation (1.7)
  we see that the system of equations (1.6)--(1.10) are mutually compatible.

\section{Invariance of the Equations (1.6)--(1.10) under the Galilean Transformation}

\hskip 2em
  In this section we prove that the equations (1.6)--(1.10) are invariant under the Galilean
  transformation.

  Let $O\bfx$ and $O'\bfx'$ be two Cartesian solid coordinate systems, where $\bfx=(x,y,z)$ and
  $\bfx'=(x',y',z')$, with $O\bfx$ fixed and $O'\bfx'$ moving in velocity $\bfv_0$ and $O'\bfx'=O\bfx$
  at time $t=0$, so that
$$
  \bfx=\bfx'+t\bfv_0, \qquad \mbox{or equivalently}, \qquad \bfx'=\bfx-t\bfv_0.
\eqno{(3.1)}
$$
  Let $\rho=\rho(\bfx,t)$ and $\rho'=\rho'(\bfx',t)$ be representations of the density of electric
  charges distributed in the space in the coordinate systems $O\bfx$ and $O'\bfx'$, respectively,
  $\bfJ=\bfJ(\bfx,t)$ and $\bfJ'=\bfJ'(\bfx',t)$ be representations of the strength of electric
  current in the coordinate systems $O\bfx$ and $O'\bfx'$, respectively, $\bfu=\bfu(\bfx,t)$ and
  $\bfu'=\bfu'(\bfx',t)$ be representations of the velocity of electric current in the coordinate
  systems $O\bfx$ and $O'\bfx'$, respectively, $\overline{\bfu}=\overline{\bfu}(t)$ and
  $\overline{\bfu}'=\overline{\bfu}'(t)$ be representations of the mean velocity of electric current
  (cf. (2.9)) in the coordinate systems $O\bfx$ and $O'\bfx'$, respectively, $\bfE=\bfE(\bfx,t)$ and
  $\bfE'=\bfE'(\bfx',t)$ be representations of the electric field strength in the coordinate systems
  $O\bfx$ and $O'\bfx'$, respectively, and $\bfB=\bfB(\bfx,t)$ and $\bfB'=\bfB'(\bfx',t)$ be
  representations of the magnetic flux density (defined as in (2.10)) in the coordinate systems
  $O\bfx$ and $O'\bfx'$, respectively. In the coordinate system $O\bfx$ we have the equations
  (1.6)--(1.10). We prove that under the variable transformation (3.1) these equations are invariant,
  i.e., in the coordinate system $O'\bfx'$ they have the following expressions:
$$
  \partial_t\bfB'+\nabla\times\bfE'=\bf0,
\eqno{(3.2)}
$$
$$
  \partial_t\bfE'-\overline{\bfu'}\times(\nabla\times\bfE')+\dot{\overline{\bfu'}}\times\bfB'
  -\nabla\times[\overline{\bfu'}\times(\bfE'+\overline{\bfu'}\times\bfB')]
  -\epsln_0^{-1}\mu_0^{-1}\nabla\times\bfB'=-\epsln_0^{-1}\bfJ',
\eqno{(3.3)}
$$
$$
  \nabla\cdot(\bfE'+\overline{\bfu'}\times\bfB')=\epsln_0^{-1}\rho',
\eqno{(3.4)}
$$
$$
  \nabla\cdot\bfB'=0,
\eqno{(3.5)}
$$
$$
  \partial_t\rho'+\nabla\cdot\bfJ'=0.
\eqno{(3.6)}
$$

  We first note that the following relations hold:
$$
  \rho'(\bfx',t)=\rho(\bfx'+t\bfv_0,t),
\eqno{(3.7)}
$$
$$
  \bfu'(\bfx',t)=\bfu(\bfx'+t\bfv_0,t)-\bfv_0, \vspace{-2mm}
\eqno{(3.8)}
$$
\begin{eqnarray*}
  \qquad\qquad\qquad\qquad\;\; \bfJ'(\bfx',t)
   &\,=\,&[\bfu(\bfx'+t\bfv_0,t)-\bfv_0]\rho(\bfx'+t\bfv_0,t)
\\
    &\,=\,&\bfJ(\bfx'+t\bfv_0,t)-\rho(\bfx'+t\bfv_0,t)\bfv_0.
  \qquad\qquad\qquad\qquad\qquad\qquad  (3.9)
\end{eqnarray*}
  To establish the relationship between $\bfE$, $\bfB$ and $\bfE'$, $\bfB'$, let us consider the
  force applied to an electric charge of quantity $Q$ moving in constant velocity $\bfv$ relative
  to the fixed coordinate system $O\bfx$. By definition of $\bfE$ and $\bfB$, the force $\bfF$ borne
  by this flowing charge $Q$ is equal to
$$
  \bfF=Q\bfE+Q\bfv\times\bfB=Q(\bfE+\bfv\times\bfB).
$$
  Since in the coordinate system $O'\bfx'$ the velocity of the charge $Q$ is $\bfv-\bfv_0$, the force
  $\bfF$ is also equal to
$$
  \bfF=Q\bfE'+Q(\bfv-\bfv_0)\times\bfB'=Q[\bfE'+(\bfv-\bfv_0)\times\bfB'].
$$
  By arbitrariness of $Q$ and $\bfv$, we get
$$
  \bfE=\bfE'-\bfv_0\times\bfB', \qquad  \bfB=\bfB',
$$
  or equivalently,
$$
  \bfE'=\bfE+\bfv_0\times\bfB, \qquad  \bfB'=\bfB,
$$
  i.e.,
$$
  \bfE'(\bfx',t)=\bfE(\bfx'+t\bfv_0,t)+\bfv_0\times\bfB(\bfx'+t\bfv_0,t),
\eqno{(3.10)}
$$
$$
  \bfB'(\bfx',t)=\bfB(\bfx'+t\bfv_0,t).
\eqno{(3.11)}
$$
  For reader's convenience we copy some formulas in the vector algebra theory and the
  mathematical field theory as follows:
$$
  \bfa\times(\bfb\times\bfc)=(\bfa\cdot\bfc)\bfb-(\bfa\cdot\bfb)\bfc,
$$
$$
  \nabla(\bfa\cdot\bfb)=(\bfb\cdot\nabla)\bfa+(\bfa\cdot\nabla)\bfb
  +\bfb\times(\nabla\times\bfa)+\bfa\times(\nabla\times\bfb).
$$
$$
  \nabla\cdot(\bfa\times\bfb)=\bfb\cdot(\nabla\times\bfa)-\bfa\cdot(\nabla\times\bfb),
$$
$$
  \nabla\cdot(f\bfa)=f\nabla\cdot\bfa+\nabla f\cdot\bfa,
$$
$$
  \nabla\times(\bfa\times\bfb)=(\bfb\cdot\nabla)\bfa-(\bfa\cdot\nabla)\bfb
  +(\nabla\cdot\bfb)\bfa-(\nabla\cdot\bfa)\bfb,
$$
$$
  \nabla\times(f\bfa)=f\nabla\times\bfa+\nabla f\times\bfa.
$$
  These formulas can be easily found from basic textbooks; see for instance the list of
  vector formulas on the first page of \cite{Jackson}.

  Now let us verify (3.2)--(3.6). First, by (3.11), (3.10), (1.6) and  (1.9) we have
\begin{eqnarray*}
  \partial_t\bfB'+\nabla\times\bfE'
  &\,=\,&[\partial_t\bfB+(\bfv_0\cdot\nabla)\bfB]+[\nabla\times\bfE
  +\nabla\times(\bfv_0\times\bfB)]
\\
  &\,=\,&(\partial_t\bfB+\nabla\times\bfE)+[(\bfv_0\cdot\nabla)\bfB
  -(\bfv_0\cdot\nabla)\bfB+(\nabla\cdot\bfB)\bfv_0]
\\
  &\,=\,&\bf0.
\end{eqnarray*}
  Hence (3.2) follows. Next, we note that clearly $\overline{\bfu'}=\overline{\bfu}-\bfv_0$, so that
  from (3.10), (3.11), (3.9), (1.6), (1.7), (1.8) and (1.9) we have
\begin{eqnarray*}
  &\,\,&\partial_t\bfE'-\overline{\bfu'}\times(\nabla\times\bfE')+\dot{\overline{\bfu'}}\times\bfB'
  -\nabla\times[\overline{\bfu'}\times(\bfE'+\overline{\bfu'}\times\bfB')]
  -\epsln_0^{-1}\mu_0^{-1}\nabla\times\bfB'+\epsln_0^{-1}\bfJ'
\\
  &\,=\,&\{[\partial_t\bfE+(\bfv_0\cdot\nabla)\bfE]+\bfv_0\times[\partial_t\bfB+(\bfv_0\cdot\nabla)\bfB]\}
  -(\overline{\bfu}-\bfv_0)\times[\nabla\times(\bfE+\bfv_0\times\bfB)]
\\
  &\,\,&+\dot{\overline{\bfu}}\times\bfB-\nabla\times\{(\overline{\bfu}-\bfv_0)\times[(\bfE+\bfv_0\times\bfB)
  +(\overline{\bfu}-\bfv_0)\times\bfB]\}
  -\epsln_0^{-1}\mu_0^{-1}\nabla\times\bfB+\epsln_0^{-1}\bfJ
\\
  &\,\,&-\epsln_0^{-1}\rho\bfv_0
\\
  &\,=\,&\{\partial_t\bfE-\overline{\bfu}\times(\nabla\times\bfE)+\dot{\overline{\bfu}}\times\bfB
  -\nabla\times[\overline{\bfu}\times(\bfE+\overline{\bfu}\times\bfB)]
  -\epsln_0^{-1}\mu_0^{-1}\nabla\times\bfB+\epsln_0^{-1}\bfJ\}
\\
  &\,\,&+(\bfv_0\cdot\nabla)\bfE+\bfv_0\times[-\nabla\times\bfE+(\bfv_0\cdot\nabla)\bfB]
  -\overline{\bfu}\times[\nabla\times(\bfv_0\times\bfB)]
\\
  &\,\,&+\bfv_0\times[\nabla\times(\bfE+\bfv_0\times\bfB)]
  +\nabla\times[\bfv_0\times(\bfE+\overline{\bfu}\times\bfB)]
  -[\nabla\cdot(\bfE+\overline{\bfu}\times\bfB)]\bfv_0
\\
  &\,=\,&\underline{(\bfv_0\cdot\nabla)\bfE}-\underline{\underline{\bfv_0\times(\nabla\times\bfE)}}
  +\underbrace{\bfv_0\times[(\bfv_0\cdot\nabla)\bfB]}
  -\overline{\bfu}\times[\nabla\times(\bfv_0\times\bfB)]
\\
  &\,\,&+\underline{\underline{\bfv_0\times(\nabla\times\bfE)}}+\underbrace{\bfv_0\times[\nabla\times(\bfv_0\times\bfB)]}
  +\underline{\nabla\times(\bfv_0\times\bfE)}+\nabla\times[\bfv_0\times(\overline{\bfu}\times\bfB)]
\\
  &\,\,&-\underline{(\nabla\cdot\bfE)\bfv_0}-[\nabla\cdot(\overline{\bfu}\times\bfB)]\bfv_0
\\
  &\,=\,&-\overline{\bfu}\times[\nabla\times(\bfv_0\times\bfB)]
  +\nabla\times[\bfv_0\times(\overline{\bfu}\times\bfB)]-[\nabla\cdot(\overline{\bfu}\times\bfB)]\bfv_0
\\
  &\,=\,&-\overline{\bfu}\times[-(\bfv_0\cdot\nabla)\bfB+(\nabla\cdot\bfB)\bfv_0]
  +\{-(\bfv_0\cdot\nabla)(\overline{\bfu}\times\bfB)+[\nabla\cdot(\overline{\bfu}\times\bfB)]\bfv_0\}
\\
  &\,\,&\;-[\nabla\cdot(\overline{\bfu}\times\bfB)]\bfv_0
\\
  &\,=\,&\bf0.
\end{eqnarray*}
  This proves (3.3).  To prove (3.4) we note that from (3.10), (3.11), (3.7) and (1.8) we have
\begin{eqnarray*}
  \nabla\cdot(\bfE'+\overline{\bfu'}\times\bfB')-\epsln_0^{-1}\rho'
  &\,=\,&\nabla\cdot[(\bfE+\bfv_0\times\bfB)+(\overline{\bfu}-\bfv_0)\times\bfB]
  -\epsln_0^{-1}\rho
\\
  &\,=\,&\nabla\cdot(\bfE+\overline{\bfu}\times\bfB)-\epsln_0^{-1}\rho=0.
\end{eqnarray*}
  Hence (3.4) follows. The equation (3.5) is an immediate consequence of (3.11) and (1.9). Finally, from
  (3.7), (3.9) and (1.10) we have
$$
  \partial_t\rho'+\nabla\cdot\bfJ'=[\partial_t\rho+(\bfv_0\cdot\nabla)\rho]
  +[\nabla\cdot\bfJ-(\bfv_0\cdot\nabla)\rho]=\partial_t\rho+\nabla\cdot\bfJ=0.
$$
  Hence (3.6) follows.

  The proof that the equations (1.6)--(1.10) are invariant under the Galilean transformation (3.1)
  is complete.

\section{An example}

\hskip 2em
  In this section we study a simple example which has been mistakenly illustrated in many textbooks
  of electromagnetism to support the idea that Galilean transformation should be replaced with Lorentz
  transformation in the study of electromagnetic issues. We shall show that if the idea in derivation
  of the equation (1.8) is employed then no paradox will appear when using Galilean transformation to
  treat electromagnetic problems.

  Consider two charges of quantities $q_1$ and $q_2$ placed at points $A$ and $B$, respectively. For
  simplicity we assume that $q_1,q_2>0$. Let $\bfr=\overrightarrow{AB}$ and $r=|\bfr|$. The Coulomb
  forces acting on $q_1$ and $q_2$ are respectively
$$
  \bfF_1=-\dfrac{q_1q_2}{4\pi\epsln_0r^3}\bfr \quad \mbox{and} \quad
  \bfF_2=\dfrac{q_1q_2}{4\pi\epsln_0r^3}\bfr.
\eqno{(4.1)}
$$
  For an observer who is motionless relative to the two charges, he only observes these forces. If the
  observer moves in velocity $\bfu$ relative to the two charges
  instead, the situation changes. For clarity of narration we only consider forces acting on the charge
  $q_2$. To this end, let us first consider electromagnetic field induced by the charge $q_1$ observed
  by the moving observer. In the eyes of the moving observer, $q_1$ induces one magnetic field and two
  electric fields. The unique magnetic field induced by $q_1$ observed by the moving observer is the
  Biot-Savart field because in the eyes of the moving observer, $q_1$ is an electric current of strength
  $-q_1\bfu\delta(P-A)$ (remember that we are only considering the value of this current at time $t=0$), 
  and the expression of this field is
$$
  \bfB_{1}(P)=-\dfrac{\mu_0 q_1}{4\pi|\overrightarrow{AP}|^3}\bfu\times\overrightarrow{AP}, \qquad
  P\in\bfK\backslash\{A\}.
$$
  The first electric field induced by $q_1$ is the Coulomb field:
$$
  \bfE_{11}(P)=\dfrac{q_1}{4\pi\epsln_0|\overrightarrow{AP}|^3}\overrightarrow{AP}, \qquad
  P\in\bfK\backslash\{A\}.
$$
  The second electric field induced by $q_1$ is that caused by the movement (in the eyes of the moving
  observer) of the magnetic field $\bfB_{1}$ (cf. (2.5)):
$$
  \bfE_{12}(P)=-(-\bfu)\times\bfB_{1}(P)=-\dfrac{\mu_0 q_1}{4\pi|\overrightarrow{AP}|^3}
  \bfu\times(\bfu\times\overrightarrow{AP}), \qquad P\in\bfK\backslash\{A\}.
$$
  Hence in the eyes of the moving observer, the resultant force on the charge $q_2$ exerted by the charge
  $q_1$ is as follows:
$$
  \bfF_2'=-q_2\bfu\times\bfB_{1}(B)+q_2\bfE_{11}(B)+q_2\bfE_{12}(B)=q_2\bfE_{11}(B)
  =\dfrac{q_1q_2}{4\pi\epsln_0r^3}\bfr=\bfF_2.
$$
  Similarly in the eyes of the moving observer, the resultant force $\bfF_1'$ on the charge $q_1$ exerted
  by the charge $q_2$ is the same as $\bfF_1$:
$$
  \bfF_1'=-\dfrac{q_1q_2}{4\pi\epsln_0r^3}\bfr=\bfF_1.
$$
  Namely, no matter the observer is motionless or is in motion, net force between the two charges $q_1$
  and $q_2$ observed by the observer is the same, and is equal to the Coulomb force between the two charges.

  In all contemporary textbooks on electrodynamics containing the above example (lectured to give rise to
  the theory of special relativity) the electric field $\bfE_{12}$ is ignored, and in the eyes  of the
  moving observer, the resultant force on the charge $q_2$ exerted by the charge $q_1$ is calculated as
  follows:
\begin{eqnarray*}
  \bfF_2''&\,=\,&q_2\bfE_{11}(B)+(-q_2\bfu)\times\bfB_{1}(B)
\\
  &\,=\,&\dfrac{q_1q_2}{4\pi\epsln_0r^3}\bfr+\dfrac{\mu_0 q_1q_2}{4\pi r^3}\bfu\times(\bfu\times\bfr)
\\
  &\,=\,&\dfrac{q_1q_2}{4\pi\epsln_0r^3}\bfr+\dfrac{\mu_0 q_1q_2}{4\pi r^3}[(\bfu\cdot\bfr)\bfu-|\bfu|^2\bfr)
\\
  &\,=\,&\dfrac{q_1q_2\bfr}{4\pi\epsln_0r^3}(1-\beta^2)+\dfrac{\mu_0 q_1q_2}{4\pi r^3}(\bfu\cdot\bfr)\bfu,
\end{eqnarray*}
  where $\beta$ is as before. In particular, if $\bfu\bot\bfr$ then $\bfF_2''=(1-\beta^2)\bfF_2$, whose
  norm is smaller than that of the Coulomb force $\bfF_2$. Namely, in the eyes of an observer who is
  moving in the direction vertical to the line connecting the two charges $q_1$ and $q_2$, interaction
  force between the two charges is smaller than that observed by a motionless observer. This is clearly
  ridiculous.

  The above problem is ``satisfactorily'' solved by introduction of the theory of special relativity into
  electromagnetics. According to this theory, when the two charges $q_1$ and $q_2$ are moving in speed
  $|\bfu|$, their charge quantities become larger than they are still, and are respectively equal to
  $q_1/\sqrt{1-\beta^2}$ and $q_2/\sqrt{1-\beta^2}$, so that
$$
  \bfF_2''(\mbox{in relativity})=\dfrac{(q_1/\sqrt{1-\beta^2})(q_2/\sqrt{1-\beta^2})}{4\pi\epsln_0r^3}\bfr(1-\beta^2)
  =\dfrac{q_1q_2}{4\pi\epsln_0r^3}\bfr=\bfF_2.
$$
  Having recognized existence of the electric field $\bfE_{12}$, we immediately see that the above inference
  is actually incorrect, and it leads to a wrong conclusion: $|\bfF_2''|$ is larger than $|\bfF_2|$. Namely,
  if we use the theory of relativity to consider this example in the right way, we actually obtain a ridiculous
  assertion: In the eyes of an observer who is moving in the direction vertical to the line connecting the two
  charges $q_1$ and $q_2$, interaction force between them is larger than that observed by a motionless observer.

  In conclusion, as far as transformation between two relatively moving frame of references is concerned,
  Galilean transformation is correct whereas Lorentz transformation  is incorrect. As a corollary, we see
  the theory of special relativity does not apply to electromagnetic phenomena in the universe
  that we humans have been living in. In a subsequent paper (see \cite{Cui}) more detailed discussion on
  this issue will be made.

\end{document}